\documentclass[a4paper, 12pt]{article}

\usepackage{graphicx}
\usepackage{amssymb}

\def\blankbox{{ ~\hfill$\rlap{$\sqcap$}\sqcup$}}

\begin{document}

\bigskip \bigskip

\centerline {\bf Some Properties of Strongly Regular Graphs}

\bigskip \bigskip

\centerline{Norman Biggs}

\bigskip

\bigskip

\centerline{Department of Mathematics}

\centerline{London School of Economics}

\centerline{Houghton Street}

\centerline{London WC2A 2AE}

\centerline{U.K.}

\centerline{n.l.biggs@lse.ac.uk}

\bigskip

\centerline{Research Report  - May 2011}

\bigskip \bigskip

\centerline{\bf Abstract}
\bigskip

An approach to the enumeration of feasible parameters for strongly regular graphs is
described,  based on the pair of structural parameters $(a,c)$ and the positive
eigenvalue $e$. The Krein bound ensures that there are only finitely many possibilities
for $c$, given $a$ and $e$, and the standard divisibility conditions can be
used to reduce the possibilities further. Many sets of feasible parameters appear to
be accidents of arithmetic, but in some cases the conditions are satisfied for
algebraic reasons.  As an example, we discuss an infinite family of feasible parameters
for which the corresponding graphs necessarily have a closed neighbourhood as a star complement for $e$.

\vfill \eject

{\bf 1. Introduction}
\medskip

A  {\it strongly regular} graph $X$ is characterized by three parameters $k$, $a$, and $c$,
according to the rules
\smallskip

{\leftskip 30pt

$\bullet$  $X$ is regular with degree $k$;

$\bullet$  any two adjacent vertices have $a$ common neighbours;

$\bullet$  any two non-adjacent vertices have $c$ common neighbours.

\par}

\smallskip

We shall discuss only cases where

$$k\ge 3, \qquad k > c \ge 1.$$

These conditions rule out some simple examples, such as the
complete bipartite graphs $K_{k,k}$, which have $c=k$.  Thus, when we say that $X$ is an
SR$(n,k,a,c)$ graph we mean that it is a non-bipartite connected graph with $n$ vertices, and
the parameters $k,a,c$ satisfy the conditions listed above.
\smallskip

Since about 1970 many lists of feasible parameters have been constructed.  Brouwer's list in the
{\it Handbook of Combinatorics} contains  all possibilities with $n \le 280$, and
his online version currently extends to $n \le 1300$ [{\bf 2}]. These lists contain some parameters
that are impossible
according to the so-called Krein bounds, and the method to be developed here will start by excluding
these cases.
\smallskip

We shall need some standard notation and theory [{\bf 9}].
If $v$ is a vertex of an SR$(n,k,a,c)$ graph $X$ we denote the subgraphs
induced by the sets of vertices at distances
$1$ and $2$ from $v$ by $X_1(v)$ and $X_2(v)$ respectively.
The sizes of $X_1(v)$, $X_2(v)$, and $X$ are given by
$$ |X_1(v)| = k, \quad |X_2(v)| = \ell =  \frac{k(k-a-1)}{c}, \quad |X| = n = 1 + k + \ell.$$
\smallskip

The complement $X^*$ of a strongly regular graph $X$ is strongly regular with $k^*  = \ell$ and
$\ell^* = k$.  We shall normally deal with a representative of the complementary pair that satisfies
$k \le \ell$.
\smallskip

The adjacency matrix $A$ of $X$ satisfies the equations

$$ AJ = kJ \qquad {\rm and} \qquad A^2 -(a-c)A - (k-c)I = cJ, $$

where $I$ is the identity matrix and
$J$ is the all-$1$ matrix. It follows that the eigenvalues of $A$ are
$k$ (with multiplicity 1) and the roots $\lambda_1, \lambda_2$
of the equation $\lambda^2 -(a-c)\lambda - (k-c) = 0$.  We shall not concern ourselves with
situation when the discriminant of this equation is irrational, because the possibilities in that case are
very limited [{\bf 9}].  So we shall assume
that there is a positive integer $s$ such that
$s^2 = (a-c)^2 + 4(k-c)$, and the eigenvalues are the integers

$$ k = \frac {s^2 - (a-c)^2}{4} + c \qquad \lambda_1 = \frac{a-c+s}{2}, \qquad \lambda_2 = \frac{a-c-s}{2}.$$

Elementary arguments show that the multiplicities $m_1, m_2$ of $\lambda_1, \lambda_2$ are given by
the formulae

$$m_1 = \frac{k}{2cs} \Big( (k+c-a-1)(s+c-a) - 2c\Big), \quad
  m_2 = \frac{k}{2cs} \Big( (k+c-a-1)(s-c+a) + 2c\Big).$$

Given that the graph $X$ exists, these formulae must represent integers.
 If we take the basic parameters to be
$a$, $c$, and $e = \lambda_1$, then $k = (e+1)c + e(e-a)$ and $s= c+2e-a$, so the conditions can be formulated
in terms of of those parameters. Further relationships linking the structural
 parameters, the eigenvalues, and the multiplicities are well known, and can be found in the standard texts
 [{\bf 7}], [{\bf 9} p.244].
A formulation suitable for our purposes will be presented in Section 3.  But first we shall look at another
condition, which effectively bounds $c$ in terms of $a$ and $e$.
\bigskip

{\bf 2. The Krein bounds and their consequences}
\medskip

Conditions asserting the non-negativity of the so-called {\it Krein parameters}
arise in the general theory of distance-regular graphs [{\bf 4}].
In the strongly regular case we have two parameters $K_1$ and $K_2$ which, after some
elementary algebra,
 can be written in terms of $k, \lambda_1, \lambda_2$ [{\bf 9}]:

$$K_1 = (k + \lambda_1)(\lambda_2 + 1)^2 - (\lambda_1 + 1)(k + \lambda_1 + 2 \lambda_1 \lambda_2).$$

$$K_2 = (k + \lambda_2)(\lambda_1 + 1)^2 - (\lambda_2 + 1)(k + \lambda_2 + 2 \lambda_1 \lambda_2).$$

Here we shall normally consider the member of a complementary pair for which the non-negativity of
$K_2$ is the effective condition.
\smallskip

The formulae given in the previous section imply that when
$\lambda_1 = e$ we have $k = e(e-a) + (e+1)c$ and $\lambda_2 = a-c-e$.
Thus we can  express $K_2$ in terms of  $a$,$c$, and $e$; in fact it is
a quadratic in $c$:

$$K_2(c) = P + Qc - ec^2$$

where  $P= (e+1)(e-a)(e^2 - e + a)$ and  $Q =  e^3 +(2a+1)e +a$.
\vfill \eject

{\bf Theorem 2.1} \quad  The Krein condition for $K_2$ implies that

$$c \le \cases{ e^2 + e + 2a &if $e\ge 3$; \cr
                e^2 + e + 3a &if $e = 1,2$. \cr}$$

In particular, when $a=0$ we have $c \le e(e+1)$.
\smallskip

{\it Proof} \quad  The Krein condition states that $K_2 \ge 0$. Since  $K_2$
is a quadratic function of $c$ with a maximum at $c = Q/2e > 0$ it follows that
$c$ must not exceed the larger root $c_0$ of the equation $K_2(c) = 0$.
\smallskip

We can estimate $c_0$ by evaluating $K_2$.  When $c=e(e+1)$ we have

$$k = e(e^2 + 3e -a+ 1) \quad {\rm and} \quad K_2 = a(e+1)(e^2 + 3e - a).$$

Since we assume that $k \ge 3$, we have $K_2(e^2 + e)  \ge 0$. There is equality
when $a=0$ so  $c_0 = e(e+1)$ in that case.
Similarly,  we find

$$ K_2(e^2 + e +a) = 2ae(e+2), $$

$$ K_2(e^2 + e+2a) = a\big( e(5-e^2) + a(1-e) \big).$$

When $a \ge 1$ the first value is positive and the
second is negative if $e \ge 3$, so $c_0$ lies strictly between $e^2 + e + a$ and $e^2 +e +2a$.
\smallskip

When $e=1$ and $e=2$ the values of $K_2(e^2 + e + 3a)$ are $2a(1-a)$ and $-6a(a+2)$ respectively, from which
the result follows.
\blankbox
\bigskip

{\bf Example 2.2} \quad In the case $a=0$ we have $c \le e^2 + e$, and the graph is an
SRNT (strongly regular, no triangles) graph.
This case was studied in earlier papers [{\bf 1}], using an approach based on the parameters $(a,c,e)$.
\medskip

{\bf Example 2.3} \quad  In the case $a=e$ we have $P=0$, and the equation $K_2(c) = 0$ has roots $0$ and

$$c_0 =  \; Q/e = \;   e^2 + 2e +2.$$

In particular, when $e= a =1$ we have $c \le 5$, and in fact graphs exist only when $c=2, 3, 5$.
The corresponding values of $n$ are $9, 15, 27$ and
the graphs are $LK_{3,3}$, $LK_6$, and the
co-Schl\"{a}fli graph [{\bf 9}].  Note that the parameter-set  $n= 63, k= 22, a= 1, c=11$  also occurs
in Brouwer's list [{\bf 2}], only to be ruled
out by the Krein bound.
\blankbox
\bigskip

{\bf 3. The integrality conditions}
\medskip

We now formulate the `standard' integrality conditions.  As in Section 2, we take
 the basic parameters to be $a, c$ and $e$.
\vfill \eject

{\bf Theorem 3.1} \quad  Suppose $X$ is an SR$(n,k,a,c)$ graph with $\lambda_1 = e$, and let

$$D =e(e+1)(e-a)(e-a-1), \quad F = (e+1)(e^2 +2e -a)(e^2 + 3e - a).$$

Then $D$ is a multiple of $c$ and $F$ is a multiple of $c+2e-a$.
\medskip

{\it Proof} \quad If $X$ exists then the `multiplicities' $m_1$ and $m_2$ must be integers.
\smallskip

Since $m_1+ m_2 = k + \ell = n-1$, the condition that $m_1 + m_2$ is an integer is equivalent to the
integrality of $\ell = |X_2(v)|$. We have

$$ \ell  = \frac{k}{c} (k-a-1), $$

and substituting $k = (e+1)c + e(e-a)$ leads to the formula
$$ \ell \; = \;  (e+1)^2 c
    +  (e+1)(2e^2 - 2ae - a - 1)
      + \frac{D}{c},$$

where $D$ is the expression given above. Thus $D$ is a multiple of $c$.
\smallskip

The formula for $m_2$ reduces to

$$m_2 \; = \;  \frac{k((k+c-a-1)e +c)}{c(c+2e-a)} \; = \; \frac{k(k-e)(e+1)}{c(c+2e-a)}.$$

Substituting for $k$, we find (after some elementary algebra)

$$m_2 \; =  \;  (e+1)^3 - \frac{Ec - De}{c(c+2e-a)},$$

where $E = (e+1)^2 (ae +3e -a)$.
\smallskip

Given that $D$ is a multiple of $c$, say $D=cc'$, the second term is an integer
if $E - ec'$ is a multiple of $c+2e-a$. If $\theta$ is the integer such that
$ E - ec'  = \theta (c+2e-a)$, then

$$ (E- \theta c) (c +2e-a) = E (2e-a) + De \hskip 130pt$$
$$ \hskip 80pt = (e+1)(e^2 + 2e -a)(e^2 + 3e -a) = F,$$

so $F$ is a multiple of $c+2e-a$, as claimed.
\blankbox
\medskip

In certain cases the conditions can be simplified.
\medskip

{\bf Corollary 3.2} If $e =a$ the sole condition is that $c+e$ divides $e^2 (e+1)^2 (e + 2)$, and if $e= a+1$
the sole condition is that $c+e+1$ divides $(e+1)^3 (e^2+ e + 1)$.
\smallskip

{\it Proof} \quad  When $e=a$ or $e=a+1$ we have $D=0$, so the first condition is trivially satisfied.
Since $c' =0$, the sole condition is that $F$ is a multiple of $c+2e-a$, which reduces to the forms stated.
\blankbox
\bigskip

{\bf 4. Enumeration of feasible parameters}
\medskip

The foregoing results  enable feasible parameters to be calculated systematically. The method is to
fix $a$ and $e$, and find those $c$ which lie in the range specified in Theorem 2.1 and satisfy the conditions

$$  c \mid D \quad {\rm and} \quad  c+2e-a \mid F.$$

{\bf Example 4.1} \quad When $a=1$ and $e=4$ we require that $1 \le c \le 22$,
 $c$ divides $120 = 2.3.4.5$, and $c+7$ divides $3105 = 3^3.5.23$.  It is easy to check that the
 only possibilities are $c= 2, 8, 20$.  The corresponding values of $(n,k)$ are
 $(243, 22)$, $(378, 52)$, and $(729,112)$.  Here we have an atypical situation, because
  graphs are known to exist for all three sets of feasible parameters [{\bf 2}].
\blankbox

\medskip

Numerical evidence suggests that, given $a$ and $e$, the number of $c$ satisfying the required conditions
is (at worst) a linear function of $e$.  However, the situation is obscured by the notoriously
complicated multiplicative structure of the integers, which means that many feasible values of $c$
are accidents of arithmetic.  More order can be imposed by regarding $e$ as an indeterminate
variable and working in the ring ${\mathbb Z}[e]$.  If we fix $a \in {\mathbb Z}$ and regard the expressions
$$D_a(e) = e(e+1)(e-a)(e-a-1),   \qquad F_a(e) = (e+1)(e^2 + 2e - a)(e^2 + 3e -a),$$

as elements of ${\mathbb Z}[e]$, then we can look for divisors $c$ of $D_a(e)$  in ${\mathbb Z}[e]$
such that $c + (2e-a)$  is a divisor of $F_a(e)$.   In the light of Theorem 2.1 we may also impose the
condition that the divisors have degree at most $2$.  When $a=0$ there are several possibilities,
 but for any positive $a$ there is only one possibility,

$$c = e(e+1), \qquad c+(2e-a) = e^2 +3e -a.$$

We shall refer to this as the {\it algebraically feasible} case, and discuss it in more detail in Section 6.
\medskip

The foregoing results also facilitate compilation of lists of strongly regular graphs in order of the number
 of vertices, $n$.
Since $n = 1 + k +\ell$, it follows from the formula for $\ell$ obtained in the proof of Theorem 3.1
 that we can write $n$ in the form

$$n  = Gc + H + \frac{D}{c},$$

where $D$ is as before and

$$G = (e+1)(e+2),$$
$$H = 2e^3 + (3-2a)e^2 - (1+4a)e - a.$$

For given $a$ and $e$, consider $n$ as a function of $c$ in the range stated in Theorem 2.1.
Clearly $n$ has only one turning point, a minimum, at the point where
$c^2 = D/G$ and $n$ takes the value $H + 2\sqrt{DG}$. Roughly speaking, the minimum occurs near the point
$c= e$, and the minimum value $n_{min}$ is about $4e^3$.
\smallskip

The general behaviour of $n$, and in particular the location of the maximum $n_{max}$, can be inferred
from a few calculations. With a few exceptions, the largest feasible value of $c$ occurs in the
algebraically feasible case $c = e(e+1)$, in which case
$n= (e^2 + 3e -a)^2$.
\medskip

These  bounds provide an effective method of
tabulating the results.  For $a=0$ (the SRNT case) and $1 \le e \le 10$ the bounds are as follows:

$$ \matrix{    e         &1 &2  &3   &4   &5   &6   &7    &8    &9    &10       \cr
            n_{min}   &4 &50 &154 &342 &638 &1066&1650 &2413 &3381 &4577   \cr
            n_{max}   &16&100&324 &784 &1600&2916&4900 &7744 &11664&16900 \cr}. $$

Suppose we wish to list all the feasible parameters for SRNT graphs with at most 1000 vertices.
According to the table, we need only carry out the calculation for $1 \le e \le 5$, since
$n_{min}(6)$ is greater than 1000. Similarly, if we list the feasible parameters for $1 \le e \le 10$,
the list will contain all possibilities with fewer than 6025 vertices, since $n_{min}(11) = 6025$.
\bigskip

{\bf 5. Star complements}
\medskip

The numerical conditions derived above are necessary, but not sufficient, for the existence of a strongly
regular graph.  For small values of the parameters they appear to be a remarkably good guide, because
graphs can be constructed for a high proportion of the feasible sets of parameters.  Whether this is true
in general remains a mystery.
\smallskip

One approach to the questions of existence and uniqueness of an SR$(n,k,a,c)$ graph is based
on the reconstruction of the graph from a hypothetical subgraph.  In this section we shall discuss the
possibility of using a {\it star complement} as the subgraph. This theory has been developed
extensively by Cvetkovi\c{c}, Rowlinson and
Simi\c{c}  [{\bf 8}].   We shall also employ arguments similar to those
used in the method of subconstituents, as described by Godsil and Royle [{\bf 9}].
\smallskip

Let $P \cup Q$ be a partition of the vertex-set of a graph $X$. The adjacency matrix of $X$ is partitioned
correspondingly in the form

$$ \pmatrix{ A_P  &B^T \cr
              B   & A_Q \cr}, $$

where $A_P$, $A_Q$, are the adjacency matrices of the subgraphs induced by $P$ and $Q$, and $B$ specifies
the edges with one vertex in $P$ and one vertex in $Q$.  In general there is no relationship between
these submatrices, but when the partition is chosen in a particular way, there is.
\smallskip

Let $X$ be a graph with $n$ vertices that has an eigenvalue $e$ with multiplicity $m$.
It is easy to show  that removing a vertex from $X$
reduces the multiplicity of $e$ by $1$, at most.  Hence the maximum cardinality of an induced subgraph
that does not
have $e$ as an eigenvalue is $n-m$. Such a subgraph is said to be a {\it star complement} for $e$
in $X$. It can be shown that star complements always exist.
For our purposes the following result [{\bf 8}] is fundamental.
\medskip

{\bf Theorem 5.1 (The Reconstruction Theorem)}  \quad  If $Q$ is a star complement for $e$ in $X$, then
the submatrices $A_P$, $A_Q$ and $B$ satisfy the equation

$$  eI - A_P  = B^T \, (eI - A_Q)^{-1} \, B. $$

(Here we adopt the
convention, strictly incorrect, of using the same notation for a subset $Q$ of the vertices
of $X$ and the subgraph induced by $Q$.
\blankbox

\medskip

In a nutshell, this equation says that the edges of $X$ having both vertices in $P$ can be reconstructed,
given $e$ and the edges with at least one vertex in $Q$.
\smallskip

Suppose we are trying to construct a strongly regular graph $X$ with parameters $(a, c, e)$.
The degree and cardinality of $X$, and the cardinality of a star complement for $e$, are thus

$$ k = (e+1)c + e(e-a), \quad n = Gc + H + c',$$

$$n-m = m_2 +1 =  \frac{k(k-e)(e+1)}{c(c+2e-a)}  + 1 ,$$

where $cc' = D$ and $D,G,H$ are given by the formulae in Sections 3 and 4.
If we can identify a set $Q$ of $n-m$ vertices that induces a subgraph which does not have
$e$ as an eigenvalue, then the construction of $X$ depends on finding a suitable matrix $B$. This defines
the edges with exactly one vertex in $P$ and, by the Reconstruction Theorem, the edges with both vertices in
$P$.   If there is a unique $B$ compatible with the given parameters, then $X$ is unique.
If there is
no such $B$, then $X$ cannot exist.
\smallskip

The fact that the reconstructed graph $X$ is strongly regular places severe restrictions on $B$.
\medskip

{\bf Lemma 5.2} \quad Let $X$ be a strongly regular graph with
parameters $(a,c,e)$ and suppose $Q$ is a star complement for $e$, with adjacency matrix $A_Q$. Then
the matrix $B$ is such that $BB^T = R$, where

$$R = cJ + e(e+c-a)I + (a-c)A_Q -  {A_Q}^2.$$

\smallskip

{\it Proof} \quad  The adjacency matrix

$$ A \;  = \; \pmatrix{ A_P  &B^T \cr
              B   & A_Q \cr}, $$

satisfies the quadratic equation given in Section 2. Taking the submatrices in the second row and column, and
substituting $k = (e+1)c + e(e-a)$, gives the result.
\blankbox
\medskip

The matrix $R$ depends only on $A_Q$,
which we assume known. The problem is to find a suitable matrix $B$ such that
$BB^T = R$, which is only possible when $R$ is positive semi-definite.  Variants of this problem are studied
in other branches of mathematics,  but our problem has some
special features.  Specifically, $B$ must be a $(0,1)$ matrix with $m$ columns.
\medskip

{\bf Example 5.3} \quad In the case $a=0$, $c=1$, $e=1$ we have $k=3$, $n=10$, and $n-m=5$.  The parameters
imply that $X$ must contain an induced $5$-cycle $Q$, for which the adjacency matrix
satisfies ${A_Q}^2 + A_Q - I = J$. Since $e=1$ is not an eigenvalue of $A_Q$, $Q$ is a star complement for $e$.
\smallskip

In this case

$$R \; = \; J + 2I -A_Q - {A_Q}^2  \; =  \; I.$$

Clearly,  a $5 \times 5$ $(0,1)$-matrix satisfying the equation $BB^T = R$  is $B=I$.
By the Reconstruction Theorem,   $I- A_P  = (I - A_Q)^{-1}$, and using the
quadratic equation for $A_Q$ it is easy to check that

$$(I-A_Q)^{-1} = A_Q + 2I - J, \quad {\rm hence} \quad A_P = J-I - A_Q.$$

This means that $P$ is a pentagram, and the familiar picture of the Petersen graph is
obtained.   A more detailed analysis leads to the conclusion that this is the only possibility, up to
isomorphism.
\blankbox
\medskip

{\bf Example 5.4} \quad In the case $a=1$, $c=3$, $e=1$ we have $k=5$, $n=15$, and $n-m=6$.
Since $K_{3,3}$ has six vertices and its eigenvalues are $3,0,-3$, we can try it as a star complement $Q$.
(In fact it can be shown that the parameters
require an induced subgraph of this form, but it is enough to
 show that the construction works on this assumption.) We have

$$A_Q = \pmatrix{ O &J \cr J &O \cr}, \quad  {A_Q}^2 =  \pmatrix{ 3J &O \cr O &3J}, $$

$$R = 3J + 3I - 2A_Q - {A_Q}^2  \; = \; \pmatrix{ 3I &J \cr J & 3I \cr}. $$

A $6 \times 9$ $(0,1)$-matrix $B$ satisfying $BB^T = R$ is

$$ B \; = \; \pmatrix{I &I &I \cr E_1 &E_2 &E_3 \cr} \quad {\rm where}$$

$$E_1 = \pmatrix{ 1 &1 &1 \cr 0 &0 &0 \cr 0 &0 &0 \cr}, \quad
E_2 = \pmatrix{ 0 &0 &0 \cr 1 &1 &1 \cr 0 &0 &0 \cr}, \quad
E_1 = \pmatrix{ 0 &0 &0 \cr 0 &0 &0 \cr 1 &1 &1 \cr}. $$

Applying the Reconstruction Theorem, we find

$$A_P = \pmatrix{ O &F &F \cr F &O &F \cr F &F &O \cr} \quad
{\rm where} \quad F = \pmatrix{0 &1 &1 \cr 1 &0 &1 \cr 1 &1 &0 \cr},$$

from which it is clear that $P =LK_{3,3}$.
\blankbox
\bigskip

{\bf 6. The algebraically feasible case}
\medskip

In some cases a suitable star complement is suggested by the structure of the graph.
For example, the {\it closed neighbourhood}
$N  = \{v\} \cup X_1(v)$, where $v$ is any vertex in a strongly regular graph with parameters $(a,c,e)$,
 consists of $k+1$ vertices, $k$ of
which have degree $a+1$ and one of which has degree $k$. The following results show that $N$ is potentially
a star complement for $e$, for infinitely many sets of feasible parameters.
\smallskip

The formula obtained in Lemma 5.2 can be reduced to a more obvious form when $Q = N$.
Suppose $w$ is any vertex in $X_1(v)$.  Then

$$(A_N)_{vv} = 0, \quad (A_N)_{vw} = 1, \quad (A_N^2)_{vv} = k, \quad (A_N^2)_{vw} = a.$$

Hence

$$R_{vv} = c + e(e+c-a) - k = 0, \quad R_{vw} = c + (a-c) - a = 0.$$

Removing row $v$ and column $v$ from $R$, $A_N$, and $A_N^2$, we get $k \times k$ matrices $R_{\flat}$,
$A_1$, and $A_2$ such that

$$R_{\flat} = cJ + e(e+c-a)I + (a-c)A_1 - A_2.$$

(Note that $A_2$ is not the same as $(A_1)^2$.)
\smallskip

The matrix $B_{\flat}$ obtained from $B$ by removing row $v$ (which must clearly be a row of $0$s) is such
that $B_{\flat} B_{\flat}^T = R_{\flat}$. So $B_{\flat}$ is simply the  $(0,1)$-matrix of size $k \times \ell$
that specifies the edges between $X_1(v)$ and the putative $X_2(v)$. The critical fact is that if $N$
is a star complement, the entire structure of $X_2(v)$ can be reconstructed from $B_{\flat}$.
\medskip

Recall that in Section 4 the parameters $(a,c,e)$ were shown to be `algebraically feasible' whenever
$c = e(e+1)$.  In that case the corresponding values of $k$ and $n$ are

$$ k = e(e^2 + 3e -a +1), \qquad n = (e^2 + 3e -a )^2.$$

\medskip

{\bf Theorem 6.1} \quad  Let $X$ be a strongly regular graph with parameters $(a,c,e)$ such that
$c = e(e+1)$. Then for all $a$ and all $e >a$ the closed neighbourhood $N$ of a vertex is a star
complement for $e$ in $X$.
\smallskip

{\it Proof} \quad  We  have $m_2 =  {k(k-e)(e+1)}/{c(c+2e-a)}$, and elementary algebra shows that
if $c =e(e+1)$ then $(k-e)(e+1) = c(c+2e-a)$.  Hence $m_2 = k$ and $n-m = m_2 + 1  =  k+1$, so
$N$ has the right size.
\smallskip

It remains to show that $e$ is not an eigenvalue of $N$. The first subconstituent $X_1$ is a regular
graph with degree $a$ and $k$ vertices.  Hence its characteristic polynomial has the form

$$P(x) = (x-a)^p R(x),$$

where $p$ is the number of components, and the zeros of $R(x)$ are such that $|x| \le a$.   The
closed neighbourhood $N$ is obtained by joining one new vertex to all vertices of $X_1$, and by
a standard result its characteristic polynomial is

$$(x-a)^{p-1} R(x) (x^2 - ax - k).$$

Hence all the eigenvalues of $N$ satisfy $|x| \le a$, except possibly for the roots of $x^2 - ax -k =0$.
However,  putting $x=e$ we have

$$e^2 - ae - k = e^2 - ae  - e(e^2 + 3e - a-1) = -e^2 (e+3) \neq 0.$$

So $e$ is not a root, and if $e>a$, then $e$ is not an eigenvalue of $N$.
\blankbox
\medskip

{\bf Example 6.2} \quad When $a=0$, $N$ is the graph $K_{1,k}$.  In the cases $e=1$ and $e=2$ the
closed neighbourhood can be used to reconstruct the graphs uniquely [{\bf 8}], resulting in
the Clebsch and Higman-Sims graphs.  However $e=3$ corresponds to a graph with $324$
vertices, and it has been shown by other means that the graph does not exist
[{\bf 11}, {\bf 15}].  When $a=1$, $N$ is a `windmill' with $k/2$ triangles, and
 in the case $e=2$
Stevanovi\c{c} and Milosevi\c{c} [{\bf 14}].
used the reconstruction method to prove that there is a unique graph.
This graph had previously been studied by several authors [{\bf 3}, {\bf 5}].
\blankbox
\medskip

In Example 6.2 we noted two graphs with $a=e-1$. When  $e=1$ we have the
Clebsch graph SR(16, 5, 0, 2), and when $e=2$ we have a graph SR(81, 20, 1, 6).
 It is natural to ask if these graphs are part of a family with algebraically
feasible parameters of the form

$$ n= (e+1)^4, \quad k = e(e^2 +2e + 2), \quad a = e-1, \quad c = e(e+1). $$

It follows from Theorem 6.1 that. for these parameters, the closed neighbourhood $N$ is a star complement for $e$.
The  matrix $B_{\flat}$ must satisfy $B_{\flat}B_{\flat}^T = R_{\flat}$ where, according to Lemma 5.2,

$$R_{\flat} = e(e+1)J + e(e^2 +e +1)I - (e^2 +1) A_1 - A_2.$$

The isomorphism type of $N$ is not determined by the parameters, so $A_1$ and $A_2$ are not generally
known. But some  progress can be made on the assumption that $N$ is a `generalized windmill', comprising
$e^2 + 2e +2$ cliques of size $e+1$ with one vertex in common. (See also [{\bf 13}, Theorem 4.5].)
 In that case $A_1$ can be written as a
block-circulant matrix with $e^2 + 2e + 2$ rows and columns of  blocks,
each block of size $e \times e$.

$$A_1 = {\rm bcirc}[J-I \quad O  \: O \: \cdots \: O].$$

$A_2$ is a  matrix of the same form:

$$A_2 = {\rm bcirc}[(e-1)J+I \quad J   \: J \: \cdots \: J].$$

It follows that

 $$R_{\flat} = {\rm bcirc}[U \quad V \: V \: \cdots  \: V], \quad U = e(e+1)^2 I, \quad V =  (e^2+e -1)J.$$

\smallskip

It is clear that in the `generalized windmill' case $X_1(v)$ consists of $e^2 +2e +2$ cliques of size $e$.
Furthermore, it follows from
results of Brouwer and Haemers [{\bf 5}] that $X_2(v)$ is also highly structured.
These authors show that the $(e+1)^4$ vertices of $X$ can be partitioned into $e+1$ subsets of size $(e+1)^3$,
in such a way that each subset induces a graph isomorphic to $(e+1)^2$ copies of $K_{e+1}$. This implies that
$X_2(v)$ contains $e(e+2)$ cliques of size $e+1$ and $e(e+1)^2$ cliques of size $e$.  The
required  $B_{\flat}$ can thus be written as a partitioned matrix conforming with this structure:

$$B_{\flat}  \; = \; [ Y \; Z],$$

where $Y$ has $e^2 + 2e +2$ rows and $e(e+2)$ columns of blocks of size $e \times (e+1)$ and
$Z$ has $e^2 + 2e +2$ rows and $e(e+1)^2$ columns of blocks of size $e \times e$.

\medskip

{\bf Example 6.3} \quad When $e=1$ the block matrices $U$ and $V$ are the singletons $[4]$ and $[1]$,
so we require $B_{\flat} B_{\flat}^T =$ circ$[4\; 1\; 1\; 1 \;1]$.
Taking

$$ Y \; = \; \pmatrix{ 00 &00 &00 \cr
                       10 &10 &10 \cr
                       10 &01 &01 \cr
                       01 &10 &01 \cr
                       01 &01 &10 \cr} \quad {\rm and} \quad
Z = \pmatrix{ 1 &1 &1 &1 \cr
              1 &0 &0 &0 \cr
              0 &1 &0 &0 \cr
              0 &0 &1 &0  \cr
              0 &0 &0 &0  \cr},$$

it is easy to check that $Y Y^T + Z Z^T  = 3I +J$, which is the required result.
Applying the Reconstruction Theorem $I - A_P = B^T (I-A_Q)^{-1} B$ it turns out
that $A_P$ is the adjacency matrix of the Petersen graph, and we get the Clebsch graph SR(16, 5, 0, 2).
\blankbox
\medskip

The structural result of Brouwer and Haemers holds for all values of $e$. However graphs are known only when
 $e+1$ is a prime power.
 The following construction is given by Brouwer and Haemers, based on a more general
method of Ivanov and Shpectorov [{\bf 10}].   Consider the field ${\mathbb F}_{q^2}$ with the
automorphism $x \mapsto {\bar x} = x^q$.  A $2 \times 2$ matrix $M = (m_{ij})$ over ${\mathbb F}_{q^2}$
is {\it Hermitian} if $M^T = {\bar M}$: explicitly this means that $m_{11}$ and $m_{22}$ are in the
ground field ${\mathbb F}_q$ and $m_{21} = {\bar m}_{12}$.  It follows that there are $q^4$ such
matrices, and we take the set ${\cal H}$ of them to be the vertices of a graph.  Note that ${\cal H}$
is a group with respect to addition, but not with respect to multiplication.
\smallskip

Let ${\cal S}$ be the subset of matrices with rank $1$, that is

$${\cal S} = \{M  \in {\cal H} \mid M \neq O \; {\rm and} \; \det M = 0\}.$$

If $M \in {\cal S}$ then there $q-1$ non-zero possibilities for $m_{11}$ and for each of them
there are $q^2$ possibilities for $m_{12}$.  The values of $m_{21}$ and $m_{22}$ are then
determined by the equations $m_{21} = {\bar m}_{12}$ and $m_{22} = {m_{11}}^{-1} m_{12} m_{21}$.
  If $m_{11} = 0$
then $m_{12} = m_{21} = 0$ and there are $q-1$ possibilities for $m_{22}$. It follows that
$| {\cal S} | = (q-1)(q^2 + 1)$.
\smallskip

If $M \in {\cal S}$ then $\alpha M \in {\cal S}$ for all $\alpha \neq 0$.  In particular,
${\cal S} = - {\cal S}$, so we can define the Cayley graph $({\cal H}, {\cal S})$, in which
$M$ and $M'$ are adjacent whenever $M - M' \in {\cal S}$.  The first subconstituent $X_1(O)$
comprises the $(q-1)(q^2 +1)$ members of ${\cal S}$. They form $q^2 +1$ cliques of size $q-1$,
each of the form $\{\alpha M \mid \alpha \neq 0\}$.
\smallskip

It can be checked that $c = q(q-1)$ and so
we have a strongly regular graph with parameters

$$  n= q^4, \quad k = (q -1) (q^2 + 1), \quad a = q -2, \quad c = q(q-1). $$
The eigenvalue $e$ is indeed $q-1$, so the parameters $(a,c,e)$ are as
postulated above.   It has been shown by various means that
the graph is unique in the case $q=3$ [{\bf 5}, {\bf 13}, {\bf 14}].
\smallskip

Other constructions of the same family are known.  For example, the
{\it affine polar graphs} are related to two-weight codes and certain
geometrical configurations  [{\bf 6}].
Here the vertices of the graph are represented by the elements $x \in {\mathbb F}_{q^4}$, and
 $x$ and $y$
are adjacent when $Q(x-y)$ is a non-zero square, for a suitable quadratic form $Q$.
Variants of this construction have been studied [{\bf 3}], leading to other
graphs with algebraically feasible parameters.
\smallskip

As far as I am aware, these constructions work only when $e$ is of the
form $q-1$ with $q$ a prime power.   The fact that the Brouwer-Haemers structure theorem
holds for all values of $e$ raises an obvious question.

\bigskip

{\bf References}
\medskip

{\bf 1.} N.L. Biggs.  Strongly regular graphs with no triangles.  September 2009, arXiv: 0911.2160v1.
Families of parameters for SRNT graphs. October 2009, arXiv: 0911.2455v1.

{\bf 2.} A.E. Brouwer.  Strongly Regular Graphs. In: {\it Handbook of Combinatorial Designs}, ed. C. Colbourn,
J. Dinitz,  CRC Press, 1996. See also: www.win.tue.nl/~aeb/srg/.

{\bf 3.} A.E. Brouwer. Some new two-weight codes and strongly regular graphs.  {\it Discrete Applied Math.}
10 (1985) 111-114.

{\bf 4.} A.E. Brouwer, A.M. Cohen, A.Neumaier.  {\it Distance-Regular Graphs}, Springer, Berlin 1989.

{\bf 5.} A.E. Brouwer, W.H. Haemers. Structure and uniqueness of the $(81,20,1,6)$ strongly regular
graph.  {\it Discrete Math.} 106/107 (1992) 77-82.

{\bf 6.} R. Calderbank, W.M. Kantor. The geometry of two-weight codes. {\it  Bull. London Math. Soc.}
18 (1986) 97-122.

{\bf 7.} P. Cameron, J. van Lint. {\it Designs, Graphs, Codes, and their Links}. Cambridge University Press,
1991.

{\bf 8.} D. Cvetkovic, P.Rowlinson, S.K. Simic. {\it Eigenspaces of Graphs}, Cambridge University
Press, 1997.

{\bf 9.} C.D. Godsil, G. Royle. {\it Algebraic Graph Theory}, Springer, New York 2001.

{\bf 10.} A.A. Ivanov, S.V. Shpectorov. Characterization of the association schemes of Hermitian forms.
{\it J. Math. So. Japan} 43 (1991) 25-48.

{\bf 11.} P. Kaski, P.R.J. Ostergard. There are exactly five biplanes with $k=11$.  {\it J. Combinat. Des.}
16 (2008) 117-127.

{\bf 12.} M. Milosevic. An example of using star complements in classifying strongly regular graphs.
{\it Filomat} 22 (2008) 53-57.

{\bf 13.}  P. Rowlinson, B. Tayfeh-Rezaie.  Star complements in regular graphs: old and new results.
{\it Linear Algebra and its Applications} 432 (2010) 2230-2242.

{\bf 14.} D. Stevanovic, M. Milosevic. A spectral proof of the uniqueness of the strongly regular graph
with parameters (81, 20, 1, 6). {\it European J. Combinatorics} 30 (2009) 957-968.

{\bf 15.}  A.L. Gavrilyuk, A.A. Makhnev.  On Krein graphs without triangles. {\it Doklady Mathematics}
72 (2005) 591-594.
\end{document}